\documentclass[12pt]{amsart}
\usepackage{diagrams,amssymb}

\diagramstyle[height=2em,width=2em,scriptlabels]
\newarrow{Epi}----{>>}
\newarrow{Mono}>--->
\newarrow{Iso}>---{>>}
\newarrow{Mapsto}|--->
\newarrow{Igual}=====
\newarrow{Dashto}{}{dash}{}{dash}>

\newcommand{\Z}{{\mathbb Z}}
\newcommand{\Q}{{\mathbb Q}}

\newcommand{\F}{{\mathbb F}}
\newcommand{\W}{{\mathcal W}}

\newcommand{\pcom}{_{p}^{\wedge}}
\newcommand{\doscom}{_{2}^{\wedge}}

\newcommand{\hocolim}{\operatornamewithlimits{hocolim}}

\newcommand{\map}{\operatorname{Map}\nolimits}

                     \newcommand{\Rep}{\operatorname{Rep}\nolimits}

\newcommand{\A}{\ifmmode{\mathcal{A}}\else${\mathcal{A}}$\fi}
\newcommand{\K}{\ifmmode{\mathcal{K}}\else${\mathcal{K}}$\fi}
\newcommand{\U}{\ifmmode{\mathcal{U}}\else${\mathcal{U}}$\fi}
\newtheorem{Thm}{Theorem}[section]
\newtheorem{Prop}[Thm]{Proposition}

\newtheorem{Lem}[Thm]{Lemma}
\newtheorem{Conj}[Thm]{Conjecture}

\theoremstyle{definition}
\newtheorem{Defi}[Thm]{Definition}
\newtheorem{Rmk}[Thm]{Remark}

\theoremstyle{remark}

\newcommand{\ord}{\operatorname{ord}\nolimits}
\newcommand{\Ind}{\operatorname{Ind}\nolimits}
\newcommand{\GL}{\operatorname{GL}\nolimits}

\newcommand{\Mor}{\operatorname{Mor}\nolimits}

\title[Symplectic groups are $N$-determined \emph{2}-compact groups]
{Symplectic groups are $N$-determined \emph{2}-compact groups}

\thanks{The first
author is partially supported by the Ministry for Education, Science and Sport of
Republic of Slovenia Research Program no.~101-509. The second author is partially
supported by the DGES-FEDER grant BFM2001-1825, and Junta de Andaluc{\'\i}a Grant
FQM-0213.}

\author{Ale\v s Vavpeti\v c} \address{Fakulteta za matematiko in fiziko\\
Univerza v Ljubljana\\ Jadranska 19\\ SI-1111 Ljubljana\\Slovenia}
\email{ales.vavpetic@FMF.Uni-Lj.Si}
\author{Antonio Viruel}
\address{Dpto de {\'A}lgebra, Geometr{\'\i}a y Topolog{\'\i}a\\
Universidad de M{\'a}laga\\ Apdo correos 59 \\ 29080 M{\'a}laga
\\Spain} \email{viruel@agt.cie.uma.es}

\date{\today}

\begin{document}
\maketitle

\begin{abstract}
We show that for $n\ge 3$ the symplectic group $Sp(n)$ is as a $2$-compact group determined up to
isomorphism by the isomorphism type of its maximal torus normalizer.
This allows us to determine the integral homotopy type of $Sp(n)$
among connected finite loop spaces with maximal torus.
\end{abstract}

\section{Introduction}

The advent of $p$-compact groups in the celebrated work of Dwyer and Wilkerson \cite{DW} is the
culmination of a research program that can be traced back to the work of Hopf and Serre on
$H$-spaces and loop spaces, and fits within the philosophy of Hilbert's Fifth Problem: which are
the non differential (here homotopy theoretical) properties that characterize compact Lie groups?

A $p$-compact group is a loop space $(X,BX,e)$, i.e.\ $e:X\simeq \Omega(BX)$ for a pointed space
$BX$, such that $H^*(X;\F_p)$ is finite and $BX$ is $p$-complete in the sense of Bousfield and Kan
\cite{BK}. As expected, examples of $p$-compact groups are given by $p$-completion of compact Lie
groups $G$ for which $\pi_0G$ is a $p$-group, since $G\pcom$ is homotopy equivalent to
$\Omega(BG\pcom)$. In this way a $p$-compact torus $T$ of rank $n$ is the $p$-completion of an
ordinary torus, hence $BT$ is the Eilenberg-MacLane space $K((\Z_p)^{\oplus n},2)$. Further
examples are given by the realization of polynomial algebras, i.e.\ loop spaces $\Omega BX$, where
$BX$ is $p$-complete and has polynomial mod-$p$ cohomology (\cite{A}, \cite{CE}, \cite{DW3},
\cite{N3}, \cite{Q}, \cite{Z}). The importance of $p$-compact groups consists of a dictionary
(reviewed in Section \ref{dict}) that translates much of the rich internal algebraic structure of
compact Lie groups to the homotopy theoretical setting of $p$-compact groups, so the challenge is
then to give homotopy theoretical proofs of classical algebraic Lie group theory results.

One of those challenges quoted above is the following: $p$-compact groups admit maximal tori, Weyl
groups and maximal torus normalizers in a way that extends the classical concepts in Lie group
theory \cite[8.13 \& 9.5]{DW}, so can we ``reproof" the Lie group theoretical
Curtis-Wiederhold-Williams' theorem \cite{CWW} in the setting of $p$-compact groups? Recall that
Curtis-Wiederhold-Williams' theorem states that two compact connected Lie groups are isomorphic if
and only if their maximal torus normalizers are isomorphic, hence we are led to the following
conjecture \cite[Conjecture 5.3]{D}

\begin{Conj}\label{conjetura} Let $X$ be a connected $p$-compact group with maximal torus $T_X$. Then $X$ is
determined up to equivalence by the loop space $NT_X$.
\end{Conj}

We shall say that a $p$-compact group $X$ is $N$-determined if $X$ verifies Conjecture
\ref{conjetura} even if the hypothesis ``connected" is dropped, i.e.\ $X$ is $N$-determined if
every $p$-compact group $Y$, with the normalizer of a maximal torus isomorphic to that of $X$, is
isomorphic to $X$.

Given an odd prime $p>2$, $p$-compact groups are known to be $N$-determined \cite{AGMV}, what
leads to the classification of $p$-compact groups for $p$ odd. But the situation is quite
different at $p=2$: there exist $2$-compact groups which are not $N$-determined. For example
$O(n)\doscom$ and $SO(n+1)\doscom$ are non isomorphic $2$-compact groups that have isomorphic
maximal torus normalizers. So at $p=2$ we cannot drop the hypothesis ``connected" in Conjecture
\ref{conjetura}.

We say that a $2$-compact group $X$ is weakly $N$-determined if every $2$-compact group $Y$, for
which there exists a homotopy equivalence $BN_X\simeq BN_Y$ between the maximal torus normalizers
of $X$ and $Y$, inducing an isomorphism $\pi_0X\cong \pi_0Y$, is isomorphic to $X$. From the
definitions it follows that an $N$-determined $2$-compact group is also weakly $N$-determined.

It has been shown that the $2$-compact groups $O(n)\doscom$, $SO(2n+1)\doscom$ and
$Spin(2n+1)\doscom$ \cite{N1} are weakly $N$-determined $2$-compact groups (and they are not
$N$-determined), and that $U(n)\doscom$ for $n\not= 2$ \cite{MN1}, $(G_2)\doscom$ \cite{V},
$(F_4)\doscom$ \cite{VV}, and $DI(4)$ \cite{N55} are $N$-determined ($U(2)\doscom$ is only weakly
$N$-determined, because the normalizer $N$ of a maximal torus of $U(2)\doscom$ is also a
$2$-compact group but $N$ is not isomorphic to $U(2)\doscom$). In this paper we prove that the
symplectic groups $Sp(n)\doscom$ are $N$-determined $2$-compact groups for $n\ge 3$.

\begin{Thm}\label{elteo}
Let $n\ge 3$ and let $X$ be a $2$-compact group with the maximal torus normalizer
$f_N\colon N\rTo X$ isomorphic to that of $Sp(n)\doscom$. Then $X$ and $Sp(n)\doscom$
are isomorphic $2$-compact groups.
\end{Thm}
\begin{proof}
First we prove that $X$ is connected in Section \ref{Connectedness} (Proposition \ref{connected}).
In Section \ref{coh} we show that mod-$2$ cohomology of $BX$ is isomorphic to that of $BSp(n)$ as
algebras over the Steerod algebra, which implies that the Quillen categories associated to $X$ and
$Sp(n)$ are isomorphic. In section \ref{stubborn} we describe the $2$-stubborn decomposition of
the group $Sp(n)$, which allows us to define a map from $BSp(n)\doscom$ to $BX$ that happens to be
an equivalence. This is done in Section \ref{map}.
\end{proof}

Notice that the hypothesis $n\ge 3$ is necessary as $Sp(1)\doscom=SU(2)\doscom$ and
$Sp(2)\doscom=Spin(5)\doscom$ are only weakly $N$-determined $2$-compact groups.

The combination of the results in \cite{AGMV} and Theorem \ref{elteo} shows that if $G$ is a
connected compact Lie group, then $BG$ is in the adic genus of $BSp(n)$ if and only if $G=Sp(n)$,
what in view of \cite{N6} characterize the integral homotopy type of $BSp(n)$ as a loop space.
Thus our final result is

\begin{Thm}
Let $L$ be a connected finite loop space with a maximal torus normalizer isomorphic to
that of $Sp(n)$. Then $BL$ is homotopy equivalent to $BSp(n)$.
\end{Thm}

{\bf Notation.} Here all spaces are assumed to have the homotopy type of a CW-complex. Completion
means Bousfield-Kan completion \cite{BK}. For a given space $X$, we write $H^*X$ for the mod-$2$
cohomology $H^*(X;\F_2)$. For a prime $p$, we write $X\pcom$ for the Bousfield-Kan $p$-completion
($(\Z_p)_\infty$-completion in the terminology of Bousfield and Kan) of the space $X$. We assume
that the reader is familiar with Lannes' theory \cite{L}.

\section{The dictionary} \label{dict}

As announced in the introduction, this section is devoted to a brief review of the dictionary
translating constructions and arguments from the algebraic theory of groups to the homotopical
setting of $p$-compact groups. The aim of the minimalist style of this section is to ease the
search of concepts by the reader who will find a more detailed exposition in the original
\cite{DW1}, or the reviews \cite{D}, \cite{M3} and \cite{N4} if needed.

Along this section $X$ and $Y$ will denote $p$-compact groups whose classifying spaces are $BX$
and $BY$ respectively. By $T$ we shall denote a $p$-compact torus, i.e. $BT\simeq K(\Z_p^n,2)$
where $n$ is the rank of $T$. Finally, we define:

\begin{itemize}
\item {\bf Homomorphisms \cite[\S 3.1]{DW1}:} A homomorphism $X\rTo^f Y$
of $p$-compact groups is a pointed map $BX\rTo^{Bf} BY$. The homomorphism $f$ is an isomorphism if
$Bf$ is a homotopy equivalence. It is a monomorphism if the homotopy fiber $Y/X$ of $Bf$ is
$\F_p$-finite or equivalently if $H^*(BX,\F_p)$ is a finitely generated module over $H^*(BY,\F_p)$
via $Bf^*$.

\item {\bf Centralizers \cite[\S 3.4]{DW1}:} For a homomorphism $Y\rTo^f X$
of $p$-compact groups, the centralizer $C_X\big(f(Y)\big)$ is defined by the equation
$BC_X\big(f(Y)\big):=\map(BY,BX)_{Bf}$.

\item {\bf Maximal tori \cite[Definition 8.9]{DW1}:} A monomorphism $T\rMono X$ of
a $p$-compact torus into a $p$-compact group $X$ is a maximal torus if $C_X(T)$ is a $p$-compact
toral group and if $C_X(T)/T$ is homotopically discrete. Every $p$ compact group admits maximal
tori \cite[Theorem 8.13]{DW1}.

\item {\bf Weyl group \cite[Definition 9.2]{DW1}:} Let $BT_X\rTo^{Bf_T} BX$
be a maximal torus of a $p$-compact group $X$. Assume that $Bf_T$ is already a fibration and treat
$\W_X$ as the space of self-maps of $BT_X$ over $BX$. Composition gives $\W_X$ the structure of an
associative topological monoid. It is shown \cite[Proposition 9.5]{DW1} that $\W_X$ is
homotopically discrete and therefore $W_X:=\pi_0\W_X$ is a (finite) group. Moreover, if $X$ is
connected, the action of $W_X$ on $BT_X$ induces a faithful representation $$W_X\rMono
\GL(H^*_{\Q_p}BT_X)\cong \GL_n(\Q^\wedge_p)$$ whose image is generated by pseudoreflections, i.e.
$W_X$ is a pseudo reflection group \cite[Theorem 9.7]{DW1}.

\item {\bf Maximal torus normalizers \cite[Definition 9.8]{DW1}:}
Let $BT_X\rTo^{Bf_T} BX$ be a maximal torus of a $p$-compact group $X$. The normalizer of $T_X$
denoted by $NT_X$, or simply by $N_X$ or $N$, is the loop space such that $BNT_X$ is the Borel
construction associated to the action of $\W_X$ on $BT_X$.
\end{itemize}

All these concepts generalize the classical algebraic definitions. In particular, if $G$ is a
compact Lie group such that $\pi_0G$ is $p$-group, $i\colon T\rTo G$ a maximal torus of $G$, $W$
the Weyl group of $G$, and $N$ is the normalizer of the maximal torus $T$, then the $p$-completion
$i\pcom\colon T\pcom\rTo G\pcom$ is a maximal torus of the $p$-compact group $G\pcom$. The Weyl
group $W$ is naturally isomorphic to the Weyl group $W_{G\pcom}$. The classifying space $BN$ of
the normalizer $N$ sits in the fibration $BT\rTo BN\rTo BW$, and a normalizer of the maximal torus
$T\pcom$ of $G\pcom$ is isomorphic to the fiberwise $p$-completion $BN^{\circ}_p$ by
\cite[Proposition 1.8]{MN1}, or \cite[Lemma 6.1]{V}.

\section{Connectedness}\label{Connectedness}

In this section we proceed with the first step in the proof of Theorem \ref{elteo} by proving the
following proposition.

\begin{Prop}\label{connected}
Let $X$ be a $2$-compact group with the normalizer of a maximal torus isomorphic to
that of $Sp(n)\doscom$, where $n\ge 3$, then $X$ is connected.
\end{Prop}

The proof of the result above requires calculating the Weyl group of some centralizer in the
connected component of $X$. This is done by means of the technics developed by Dwyer and Wilkerson
in \cite{DW1} that we recall now.

An extended $p$-discrete torus $P$ is an extension of a $p$-discrete torus $(\Z/p^\infty)^n$ by a
finite group. A discrete approximation for an extended $p$-compact torus $P$ is a homomorphism
$f\colon \breve P\rTo P$, where $\breve P$ is an extended $p$-discrete torus and $f$ induces an
isomorphism $\breve P/\breve{P_0}\rTo \pi_0P$ and an isomorphism $H^*BP_0\rTo H^*B\breve{P_0}$.
Every extended $p$-compact torus has a discrete approximation \cite[Proposition 3.13]{DW1}.

\begin{Defi} \cite[Definition 7.3]{DW1}
If $s\in W$ is a pseudoreflection of order $\ord(s)$, then
\begin{enumerate}
\item the {\it fixed point set} $F(s)$ of $s$ is the fixed point set of the action of $x$
on $\breve T$ by conjugation, where $x\in \breve N(T)$ is an element which projects on $s$ by
natural projection $\breve N(T)\rTo W$,
\item the {\it singular hyperplane} $H(s)$ of $s$ is the maximal divisible subgroup of
$F(s)$ (so $H(s)\cong (\Z/p^\infty)^{r-1})$,
\item the {\it singular coset} $K(s)$ of $s$ is the subset of $\breve T$ given by elements
of the form $x^{\ord(s)}$, as $x$ runs through elements of $\breve N(T)$, which project to $s$ in
$W$, and
\item the {\it singular set} $\sigma(s)$ of $s$ is the union $\sigma(s)=H(s)\cup K(s)$.
\end{enumerate}
\end{Defi}

Notice that there are inclusions $H(s)\subset \sigma(s)\subset F(s)$ \cite[Remark 7.7]{DW1}.

Let $A\subset \breve T$ be a subgroup. Let $W_X(A)$ denote the Weyl group of $C_X(A)$, and
$W_X(A)_1$ the Weyl group of the unit component $C_X(A)_0$ of $C_X(A)$. There are inclusions
$W_X(A)_1\subset W_X(A)\subset W$, where the last follows from \cite[\S 4]{DW1}. The next theorem
tells how to calculate $W_X(A)$ and $W_X(A)_1$.

\begin{Thm} \cite[Theorem 7.6]{DW1} \label{CalcCenter}
Let $X$ be a connected $p$-compact group with maximal torus $T$ and Weyl group $W$. Suppose that
$A\subset \breve T$ is a subgroup. Then
\begin{enumerate}
\item $W_X(A)$ is the subgroup of $W$ consisting of the elements which,
under the conjugation action of $W$ on $\breve T$, pointwise fix the subgroup $A$, and
\item $W_X(A)_1$ is the subgroup of $W_X(A)$ generated by those
elements $s\in W_X(A)$ such that $s\in W$ is a reflection and $A\subset \sigma(s)$.
\end{enumerate}
\end{Thm}

Now we have all the ingredients needed for the proof of Proposition \ref{connected}

\begin{proof}[Proof of Proposition \ref{connected}]
Let $X_0$ be the unit component of $X$, and let $W_{X_0}$ be the Weyl group of $X_0$. Then
$W_{X_0}$ is a normal subgroup of $W_X$ of index a power of $2$ by \cite[Proposition 3.8]{MN}. The
minimal normal subgroup of $W_X$ of 2 power index, usually denoted by $O^2(W_X)$ in the
literature, equals $(\Z/2)^{n-1}\rtimes A_n$, i.e. the sequence
$$ (\Z/2)^{n-1}\rtimes A_n=O^2(W_X)\rMono (\Z/2)^n\rtimes \Sigma_n\rEpi^\pi (\Z/2)^2, $$ where $A_n$ is the
alternating group, is exact. The group $(\Z/2)^2$ has $5$ subgroups: the trivial subgroup $1$, the
first and the second factor $Z_1$ and $Z_2$, the diagonal $D$, and the whole group $(\Z/2)^2$.
Hence, there are $5$ normal subgroups of $W_X$ of index a power of $2$:
\begin{enumerate}
\item $\pi^{-1}(1)=(\Z/2)^{n-1}\rtimes A_n$,
\item $\pi^{-1}(Z_1)=(\Z/2)^n\rtimes A_n$,
\item $\pi^{-1}(Z_2)=(\Z/2)^{n-1}\rtimes \Sigma_n$,
\item $\pi^{-1}(D)$, and
\item $\pi^{-1}((\Z/2)^2)=W_X$.
\end{enumerate}
Because $W_{X_0}$ is the Weyl group of a connected $2$-compact group, $W_{X_0}$ is a
pseudoreflection group. According to the Clark-Ewing list \cite{CE}, only the cases (3) and (5)
are pseudoreflection groups (note that $n\ge 3$). We complete the proof showing that the case
$W_{X_0}=(\Z/2)^{n-1}\rtimes \Sigma_n$ is not possible.

Suppose that $X$ is non conneted, and let $X_0$ be the unit component. By the arguments above
$W_{X_0}$ is $(\Z/2)^{n-1}\rtimes \Sigma_n$. Let $V$ be the subgroup of the maximal torus $T$ of
$X$ (and also of $X_0$) generated by elements $(-1,-1,1,\ldots)$, $(1,1,-1,-1,1,\ldots)$, and so
on. Then $V$ is an elementary abelian $2$-group of rank $m=[\frac n 2]$. Let us write $n=2m+r$
where $r$ is $0$ or $1$, and let $C$ denote the centralizer $C_{X_0}(V)$. By Theorem
\ref{CalcCenter}(1), we get $$ W_C=\{s\in W_0\mid s|_V=id_V\}=(\Z/2)^{n-1}\rtimes (\Z/2)^m, $$
where the subgroup $(\Z/2)^m\subset \Sigma_n$ is generated by the transpositions $\tau_{2i-1,2i}$
for $i=1,\ldots m$. Let $C_0$ be the unit component of $C$. By Theorem \ref{CalcCenter}(2), the
Weyl group of $C_0$ is
$$ W_{C_0}=\langle s\in W_C\mid s\hbox{ is a reflection and } V\subset\sigma(s)\rangle.
$$ An element $s\in W_C=(\Z/2)^{n-1}\rtimes (\Z/2)^m$ is a reflection if and only if $s$ equals
$((1,\ldots,1),\tau_{2i-1,2i})$ or $((1,\ldots,1,-1,-1,1,\ldots,1),\tau_{2i-1,2i})$ for some $i$,
where both entries ``$-1$" are in the $(2i-1)^{\hbox{th}}$ and $(2i)^{\hbox{th}}$ positions. We
analyze both cases:
\begin{itemize}
\item[-] If $s=((1,\ldots,1),\tau_{2i-1,2i})$, then $F(s)=\{(x_1,\ldots,x_n)\in
(\Z/2^\infty)^n\mid x_{2i-1}=x_{2i}\}$ and $H(s)=F(s)$. Therefore $\sigma(s)=F(s)$.

\item[-] If $s=((1,\ldots,1,-1,-1,1,\ldots,1),\tau_{2i-1,2i})$, then $F(s)=\{(x_1,\ldots,x_n)\in
(\Z/2^\infty)^n\mid x_{2i-1}=x_{2i}^{-1}, i=1,\ldots,m\}$. Hence also in this case
$H(s)=F(s)=\sigma(s)$.
\end{itemize}

Since $(-1)^{-1}=-1\in \Z/2^\infty$, the group $V$ is a subgroup of $\sigma(s)$ in both cases, and
by Theorem \ref{CalcCenter}, we get $$ W_{C_0}=\langle s\in W_C\mid s\hbox{ is a
reflection}\rangle=((\Z/2)^2)^m=(\Z/2)^{2m}. $$ Hence the normalizer of a maximal torus of $C_0$
has the form $M^m$, where $M$ is the subgroup of the normalizer of the maximal torus of
$Sp(2)\doscom$ corresponding to the subgroup $\langle ((1,1),\tau_{1,2}),((-1,-1),\tau_{1,2})
\rangle < (\Z/2)^2\rtimes \Z/2=W_{Sp(2)\doscom}$. By \cite[Theorem 6.1]{DW5} and \cite[Theorem
0.5B (5)]{DW2}, the group $2$-compact group $C_0$ splits into a product $C_0\cong
X_1\times\cdots\times X_m$, where each $X_i$ is isomorphic to $(SU(2)^2/E_i)\doscom$ for some
subgroup $E_i<Z(SU(2)^2)=(\Z/2)^2$, and $M$ is isomorphic to the maximal torus normalizer of
$X_i$. Among the $5$ possibilities for each $X_i$:
\begin{enumerate}
\item $(SU(2)\times SU(2))\doscom=Spin(4)\doscom$,
\item $(SU(2)/(\Z/2)\times SU(2))\doscom\cong(SO(3)\times SU(2))\doscom$,
\item $(SU(2)\times SU(2)/(\Z/2))\doscom\cong(SU(2)\times SO(3))\doscom$,
\item $(SU(2)\times_{\Z/2} SU(2))\doscom\cong SO(4)\doscom$, and
\item $((SU(2)\times SU(2))/(\Z/2)^2)\doscom\cong (SO(3)\times SO(3))\doscom$,
\end{enumerate}
only $SO(4)$ produces a pseudoreflection group which is equivalent to that given by $M$. But while
the maximal torus normalizer of $SO(4)$ is an split extension $T:(\Z/2\times\Z/2)$, $M$ is not.
Therefore there is no $2$-compact group $X_i$ whose maximal torus normalizer is $M$, what
contradicts our initial assumption of $X$ being non connected.
\end{proof}

\section{Mod-$2$ cohomology of the $2$-compact group $X$} \label{coh}

In this section we calculate the mod-$2$ cohomology of a $2$-compact group $X$ whose maximal torus
normalizer is isomorphic to that of $Sp(n)\doscom$. This is done under the induction hypothesis
that $Sp(m)\doscom$ is $N$-determined for $2<m<n$. Notice that we already know that $Sp(1)$ and
$Sp(2)$ are weakly $N$-determined.

First we need some information about the centralizers of elementary abelian subgroups in $Sp(n)$.
It is well known that theses centralizers are isomorphic to products $Sp(n_1)\times\cdots\times
Sp(n_k)$, where $n_1+\cdots +n_k=n$. The next lemma shows that those centralizers are
$N$-determined if each factor is so.

\begin{Lem} \label{indukcija}
Let $Sp(m)\doscom$ be an $N$-determined $2$-compact group for $2<m\le n$. Then the product
$Sp(n_1)\doscom\times\ldots\times Sp(n_k)\doscom$ is weakly $N$-determined if all $n_i\le n$ and
is $N$-determined if $2<n_i\le n$.
\end{Lem}
\begin{proof}
Let $Y$ be a $2$-compact group with maximal torus normalizer $N_Y$ isomorphic to that of
$Sp(n_1)\doscom\times\ldots\times Sp(n_k)\doscom$. If at least one $n_i$ is $1$ or $2$, assume
that $Y$ is connected. Since $N_Y$ is a product $N_1\times\ldots\times N_k$, where $N_i$ is the
normalizer of a maximal torus of $Sp(n_i)\doscom$, the space $Y$ is by \cite[Theorem 6.1]{DW5} isomorphic
to a product $Y_1\times\ldots\times Y_k$, where $N_i$ is the normalizer of a maximal torus of
$Y_i$. Hence $Y_i$ is isomorphic to $Sp(n_i)\doscom$ and therefore $Y$ is isomorphic to
$Sp(n_1)\doscom\times\ldots\times Sp(n_k)\doscom$. So $Sp(n_1)\doscom\times\ldots\times
Sp(n_k)\doscom$ is weakly $N$-determined if at least one $n_i$ is $1$ or $2$, otherwise $Y$ is
$N$-determined
\end{proof}

As $X$ and $Sp(n)\doscom$ ``share" the same maximal torus normalizer $N$, they both ``share" the
same maximal torus $T$. Let $E_T<T$ be the maximal toral elementary abelian $2$-group in both $X$
and $Sp(n)\doscom$. Call $f_{E_T}$ the monomorphism $E_T\rMono X$. Next lemma shows that $E_T$
is in fact the maximal elementary abelian subgroup of $X$ (up to conjugation).

\begin{Lem} \label{one_rank_n}
Let $g\colon E\rMono X$ be an elementary abelian subgroup of $X$. Then $g$ factors through
$f_{E_T}$.
\end{Lem}
\begin{proof}
If $g\colon E\rMono X$ is central, then by \cite[Lemma 4.1]{MN} or \cite[Theorem 1.2]{DW1} $g$
factors though $f_{E_T}$ (recall that $X$ is connected by Proposition \ref{connected}).

Now assume that $g\colon E\rMono X$ is not central, thus there exists a subgroup $V<E$ of rank $1$
which is noncentral. By \cite[Proof of Theorem 1.3]{M} there exists $\tilde{g}\colon E\rMono N$
such that $Bg\simeq f_N B\tilde{g}$, the centralizer $C_N(\tilde{g})$ is the maximal torus 
normalizer of $C_X(g)$, and $\tilde{g}\vert_V$ factors through $f_{E_T}$. Because $V$ is a 
toral subgroup, the centralizer $C_N(V)$ is the maximal torus normalizer of both 
$C_{Sp(n)\doscom}(V)$ and $C_X(V)$ by \cite[Theorem 1.3]{M}. So the calculation of 
$W_X(V)$ and $W_X(V)_1$ by means of Theorem \ref{CalcCenter} equals the calculation of 
$W_{Sp(n)\doscom}(V)$ and $W_{Sp(n)\doscom}(V)_1$ what implies that $C_X(V)$ is connected and since 
by induction, the centralizer $C_{Sp(n)\doscom}(V)= BSp(m)\doscom \times BSp(n-m)\doscom$, 
$m>0$, is weakly $N$-determined (Lemma \ref{indukcija}), then $C_X(V)$ is
isomorphic to $C_{Sp(n)\doscom}(V)$.

The map $g\colon E\rTo X$ has a lift to a map $g'\colon E\rTo C_V(X)\cong Sp(m)\doscom \times
Sp(n-m)\doscom$. Up to conjugacy every elementary abelian subgroup of $Sp(m) \times Sp(n-m)$ is
toral. Hence $g$ is toral, i.e.\ factors through $f_{E_T}$.
\end{proof}

We can calculate the centralizer of $E_T$ in $X$:

\begin{Lem}\label{cent-Et}
The centralizer $C_X(E_T)$ is isomorphic to the $2$-compact group $(Sp(1)^n)\doscom$.
\end{Lem}
\begin{proof}
As $E_T$ is toral, the centralizer $C_N(E_T)$ is the maximal torus normalizer of both
$C_{Sp(n)\doscom}(E_T)$ and $C_X(E_T)$ \cite[Proposition 3.4(3)]{M1}. So the calculation of
$W_X(E_T)$ and $W_X(E_T)_1$ by means of Theorem \ref{CalcCenter} equals the calculation of
$W_{Sp(n)\doscom}(E_T)$ and $W_{Sp(n)\doscom}(E_T)_1$ what implies that $C_X(E_T)$ is connected. Since
$C_{Sp(n)\doscom}(E_T)=(C_{Sp(n)\doscom}(E_T))\doscom=(Sp(1)^n)\doscom$ is weakly $N$-determined, 
the centralizer $C_X(E_T)$ is isomorphic to $C_{Sp(n)\doscom}(E_T)$ by Lemma \ref{indukcija}. Finally
$C_{Sp(n)\doscom}(E_T)\cong (Sp(1)^n)\doscom$.
\end{proof}

The action of $\Sigma_n<W_{Sp(n)}=W_X$ on $BE_T$ induces an action of $\Sigma_n$ on $BC_X(E_T)=
\map(BE_T,BX)_{Bf_{E_T}}\cong (Sp(1)\doscom)^n$ that permutes the copies $Sp(1)\doscom$. Define
$BY=BC_X(E_T)\times_{\Sigma_n} E\Sigma_n$ and consider the diagram
\begin{diagram}
(BSp(1)^n)\doscom & \rTo & BY &\rTo & B\Sigma_n\\ \uTo&&\uTo&&\dIgual \\ (BT)\doscom
&\rTo &\map(BT,BX)_{Bf_T}\times_{\Sigma_n} E\Sigma_n & \rTo & B\Sigma_n\\
\dIgual&&\dTo&&\dTo \\ (BT)\doscom & \rTo & \map(BT,BX)_{Bf_T}\times_{W_{Sp(n)}}
EW_{Sp(n)} & \rTo & BW_{Sp(n)}
\end{diagram}
where all rows are fibrations. The space $\map(BT,BX)_{Bf_T}\times_{W_{Sp(n)}} EW_{Sp(n)}$ is the
normalizer of the maximal torus $T$, so it is isomorphic to $B(T\rtimes W_{Sp(n)})$.
Therefore the space $\map(BT,BX)_{Bf_T}\times_{\Sigma_n} E\Sigma_n$ is isomorphic to 
$B(T\rtimes \Sigma_n)$. Which means that the middle row has a section, and hence also the
top row has a section. It follows that $BY$ is homotopic to $B((Sp(1)\doscom)^n\rtimes \Sigma_n)$.

\begin{Prop}\label{detect-elab}
The cohomology $H^*BX$ is detected by elementary abelian $2$-subgroups.
\end{Prop}
\begin{proof}
The cohomology $H^*BSp(1)^n$ is detected by elementary abelian $2$-subgroups, hence by 
\cite{GLZ}, $H^*BY$ is detected by elementary abelian subgroups. The normalizer $Bf_N$ 
factors trough the map $Bf_Y$.
Because $Bf_N^*$ is a monomorphism, $Bf_Y^*$ is a monomorphism. Hence $H^*BX$ is detected by
elementary abelian $2$-subgroups.
\end{proof}

We can now identify the algebra $H^*BX$:

\begin{Prop}
The cohomology $H^*BX$ is isomorphic to $H^*BSp(n)$ as an algebra over the mod-$2$ Steenrod algebra.
\end{Prop}
\begin{proof}
By Proposition \ref{detect-elab}, the cohomology $H^*BX$ is detected by elementary 
abelian $2$-subgroups, and by Lemma
\ref{one_rank_n}, every elementary abelian subgroup of $X$ factors through $E_T$. Therefore $H^*BX$
injects into $H^*BE_T$ and therefore into $H^*BC_X(E_T)$. If we take trivial action of $\Sigma_n$ on
$X$, the inclusion $C_X(E_T) \rTo X$ is a $\Sigma_n$-equivariant map. Hence the cohomology 
$H^*BX$ is a subalgebra of $(H^*BSp(1)^n)^{\Sigma_n}=H^*BSp(n)$. 
But $H^*(BX;\Q)=H^*(BT;\Q)^{W_X}=\Q[x_4,\ldots,x_{4n}]$,
hence the Bockstein spectral sequence associated to $H^*BX\subset H^*Sp(n)= \F_2[x_4,\ldots,
x_{4n}]$ converges to $\F_2[x_4,\ldots,x_{4n}]$, and therefore $H^*BX\cong H^*BSp(n)$.
\end{proof}

Recall that the Quillen category $Q_p(G)$ of a group $G$ at a prime $p$ is the category with
objects $(V,\alpha)$, where $V$ is a nontrivial elementary abelian $p$-group and $\alpha\colon
V\rTo G$ is a monomorphism, and $\Mor_{Q_p(G)}((V,\alpha),(V',\alpha'))$ is the set of group
morphisms $f\colon V\rTo V'$ such that $\alpha=\alpha'\circ f$. By Lannes' theory (\cite{L}) and
Dwyer-Zabrodsky theorem (\cite{DZ}, \cite{N5}), the set of monomorphisms $\alpha\colon V\rTo G$ is
in bijections with the set of morphisms $B\alpha^*\colon H^*BG\rTo H^*BV$ of unstable algebras
over the Steenrod algebra $\mathcal{A}_p$ such that $H^*BV$ is a finitely generated module over
$B\alpha^*(H^*BG)$. Hence, there is an equivalent description of the Quillen category which can be
use also for $p$-compact groups (\cite[\S 2]{DW6}). If $X$ is a $p$-compact group, then $Q_p(X)$ is the
category with objects $(V,\alpha)$, where $V$ is a nontrivial elementary abelian $p$-group and
$\alpha\colon H^*BX\rTo H^*BV$ is a monomorphism of unstable algebras over the Steenrod algebra
$\mathcal{A}_p$ such that $H^*BV$ is a finitely generated module over $\alpha^*(H^*BX)$, and
$\Mor_{Q_p(G)}((V,\alpha),(V',\alpha'))$ is a set of group morphisms $f\colon V\rTo V'$ such that
$\alpha=Bf^*\alpha'$. If $X$ is a compact Lie group than the definitions 
agree (\cite[Proposition 2.2]{DW6}).

By the definition of the Quillen category and propositions \ref{detect-elab}
we get the following proposition.

\begin{Prop} \label{Quillen}
The categories $Q_2(Sp(n))$ and $Q_2(X)$ are isomorphic.
\end{Prop}

\section{2-stubborn decomposition of $Sp(n)$} \label{stubborn}

A $2$-stubborn subgroup of a Lie group $G$ is a $2$-toral group $P$ such that $N_G(P)/P$ is a
finite group which has no nontrivial normal $2$-subgroup. Let $\mathcal{R}_2(Sp(n))$ be the
$2$-stubborn category of $Sp(n)$, which is the full subcategory of the orbit category of $Sp(n)$
with objects $Sp(n)/P$, where $P\subset Sp(n)$ is a $2$-stubborn subgroup. The natural map
$$ \hocolim_{Sp(n)/P\in \mathcal{R}_2(Sp(n))} ESp(n)/P \rTo BSp(n) $$ induces an isomorphism of
homology with $\Z_{(2)}$-coefficients \cite[Theorem 4]{JMO}. Hence, defining a map $f\colon
Sp(n)\doscom \rTo X$ is equivalent to defining a family of compatible maps $\{f_P\colon
ESp(n)/P\simeq BP\rTo X\mid Sp(n)/P\in ob (\mathcal{R}_2(Sp(n)))\}$.

We first recall the $2$-stubborn subgroups of $Sp(n)$ calculated in \cite{O}. Let the permutations
$\sigma_0,\ldots,\sigma_{k-1}$ in $\Sigma_{2^k}$ be defined by $$ \sigma_r(s)= \begin{cases}
s+2^r;& s\equiv 1,\ldots, 2^r \mod 2^{r+1}\\ s-2^r;& s\equiv 2^r+1,\ldots,2^{r+1} \mod 2^{r+1}
\end{cases}
$$ Let $A_0,\ldots,A_{k-1}\in Sp(2^k)$ be diagonal matrices with $$
(A_r)_{ss}=(-1)^{[\frac{s-1}{2^r}]}, $$ where $[-]$ denotes greatest interger, and
$B_0,\ldots,B_{k-1}$ be the permutation matrices for the $\sigma_0,\ldots,\sigma_{k-1}$.

\begin{Defi} \label{StubbDef}
For every $k\ge 0$, the subgroups $E_{2^k}\subset\Sigma_{2^k}$ and $\Gamma_{2^k},
\overline\Gamma_{2^k}\subset Sp(2^k)$ are defined by $$
\begin{array}{ll}
E_{2^k}&=\langle\sigma_0,\ldots,\sigma_{k-1}\rangle\cong(\Z/2)^k,\\
\Gamma_{2^k}&=\langle uI,A_r,B_r\mid u\in Q(8),0\le r< k\rangle,\\
\overline\Gamma_{2^k}&=\langle uI,A_r,B_r\mid u\in S^1(j),0\le r< k\rangle,\\
\end{array}
$$
where $Q(8)=\{\pm 1,\pm i,\pm j,\pm k\}$ is the quaternion group and 
$S^1(j)=\{a+bi,aj+bk\mid a^2+b^2=1\}$ is the normalizer of the maximal torus in $Sp(1)=S^3$.
\end{Defi}

\begin{Rmk}\label{Oli}
Let $P$ be $\Gamma_{2^k}$ or $\overline\Gamma_{2^k}$ and $P_D$ subgroup of all diagonal matrices in $P$.
Then $P_D$ is $Q(8)\times E_{2^k}$ or $S^1(j)\times E_{2^k}$ and the exstension $P_D\rTo P\rTo (\Z/2)^k$
splits.
\end{Rmk}

\begin{Thm}\cite[Theorem 3]{O}
\begin{enumerate}
\item A $2$-stubborn group $P< Sp(n)$ is irreducebile if it is conjugate to either
$$ P=\Gamma_{2^k}\wr E_{2^{r_1}}\wr\cdots\wr E_{2^{r_s}}\text{ or }\, P=\overline\Gamma_{2^k}\wr E_{2^{r_1}}\wr\cdots\wr
E_{2^{r_s}}$$ where $n=2^{k+r_1+\cdots+r_s}$.

\item A group $P<Sp(n)$ is a $2$-stubborn group if it is conjugate
to $P_1\times\cdots\times P_s$, where $P_i$ is an irreducible $2$-stubborn subgroup of
$Sp(n_i)$ and $n=n_1+\cdots +n_s$.
\end{enumerate}
\end{Thm}

Let $\widetilde{\mathcal R}_2(Sp(n))$ be the full subcategory of ${\mathcal R}_2(Sp(n))$ with
objects $Sp(n)/P$, where $P$ is one of the representative $2$-stubborn groups from the previous
theorem. The category $\widetilde{\mathcal R}_2(Sp(n))$ is equivalent to ${\mathcal R}_2(Sp(n))$,
so the natural map $$ \hocolim_{Sp(n)/P\in \widetilde{\mathcal{R}}_2(Sp(n))} ESp(n)/P \rTo BSp(n)
$$ is also a homotopy equivalence up to $2$-completion.

\begin{Prop} \label{5.5}
Let $Sp(n)/P\in \widetilde{\mathcal R}_2(Sp(n))$ and define $P_D=P\cap Sp(1)^n$ and 
$P_T=P\cap T_{Sp(n)}$. Then
\begin{enumerate}
\item $C_{Sp(n)}(P_T)=T_{Sp(n)}$ and $C_{Sp(n)}(P_D)=(\Z/2)^n$,
\item for any extension $\alpha\colon P\rTo Sp(n)$ of $i\colon P_T\rTo Sp(n)$, we have 
$C_{Sp(n)}(\alpha(P))= Z(P)$ and
\item the canonical map
$$
\pi_0(\map(BP,BSp(n)\doscom)_{B\alpha|_{BP_T}=Bi_{P_T}})\rTo Hom(H^*BSp(n),H^*BP)
$$
is an injection.
\end{enumerate}
\end{Prop}

\begin{Rmk}
By $\map(BP,BSp(n)\doscom)_{B\alpha|_{BP_T}=Bi_{P_T}}$ we denote the components of the mapping space 
$\map(BP,BSp(n)\doscom)$ given by maps $B\alpha\colon BP\rTo Sp(n)\doscom$,
such that $B\alpha|_{BP_T}\simeq Bi_{P_T}$.
\end{Rmk}

\begin{proof}
Part (1) is obvious for $P=\Gamma_{2^k}$ and $P=\overline\Gamma_{2^k}$.
If $P=Q\wr E_{2^r}$, where $Q$ is irreducible $2$-stubborn subgroup of $Sp(2^{k-r})$,
then $C_{Sp(2^k)}(P_T)=C_{Sp(2^{k-r})}(Q_T)^{2^r}$ which is, by induction,
$(T_{Sp(2^{k-r})})^{2^r}=T_{Sp(2^k)}$.
If $P=P_1\times\ldots\times P_s$ is a product of irreducible $2$-stubborn groups,
then $C_{Sp(n)}(P_T)=C_{Sp(n_1)}((P_1)_T)\times\ldots\times C_{Sp(n_s)}((P_s)_T)=
T_{Sp(n_1)}\times\ldots\times T_{Sp(n_s)}=T_{Sp(n)}$.
Analogouslly we prove that $C_{Sp(n)}(P_D)=(\Z/2)^n$.

Every map $g\in \map(BP,BSp(n)\doscom)$ is homotopic to $B\alpha\colon BP\rTo
BSp(n)\doscom$, where $\alpha$ is the $2$-completion of a homomorphism of groups
$P\rTo Sp(n)$ (\cite{DZ}, \cite{N5}).

Let $P$ be an irreducible $2$-stubborn subgroup of $Sp(2^k)$ and let $\alpha\colon P\rTo Sp(2^k)$
be a homomorphism such that $\alpha|_{P_T}=i_{P_T}$. The extensions 
$B\alpha|_{BP_D}\colon BP_D\rTo BSp(2^k)\doscom$ of $Bi_{P_T}$ are classified by
obstruction groups $H^m(P_D/P_T;\pi_m(\map(BP_T,BSp(2^k)\doscom)_{Bi_{P_T}}))$. By \cite{JMO},
$\map(BP_T,BSp(2^k)\doscom)_{Bi_{P_T}}$ is homotopy equivalent to $BC_{Sp(2^k)}(P_T)\doscom$,
which is isomorphic to $(BS^1)^{2^k}$, by part (1). Then 
\begin{equation*}
H^m(P_D/P_T;\pi_m(\map(BP_T,BSp(2^k)\doscom)_{Bi_{P_T}}))= H^m(P_D/P_T;\pi_m (BS^1)^{2^k})
\end{equation*}
and the only possible nontrivial group is for $m=1$. And
\begin{equation*}
H^1(P_D/P_T;\pi_m(\map(BP_T,BSp(2^k)\doscom)_{Bi_{P_T}}))=H^1(\Z/2;(\Z\doscom)^{2^k}),
\end{equation*}
where the group $\Z/2$ acts on $(\Z\doscom)^{2^k}$ by reflection on each component;
this action comes from the action of the Weyl group $\Z/2$ of the group $Sp(1)$
on the maximal torus $S^1$. By Shapiro's lemma \cite[III, Proposition 6.2]{B},
the group $H^1(\Z/2;(\Z\doscom)^{2^k})$ is trivial, so all obstruction groups vanish. 
Hence if $B\alpha|_{BP_T}=Bi_{P_T}$ then  $B\alpha|_{BP_D}=Bi_{P_D}$.

First we will prove part $(2)$ and $(3)$ for the case $P$ is $\Gamma_{2^k}$ or $\overline\Gamma_{2^k}$.
Let $\alpha\colon P\rTo Sp(2^k)$ be a homomorphism such that $\alpha|_{P_T}=i_{P_T}$.
Then by above paragraph $B\alpha|_{BP_D}$ is homotopic to $Bi_{P_D}$.
The extensions $B\alpha\colon BP\rTo BSp(2^k)\doscom$ of $Bi_{P_D}$ are classified by
obstruction groups $H^m(P/P_D;\pi_m(\map(BP_D,BSp(2^k)\doscom)_{Bi_{P_D}}))$. By \cite{JMO},
the mapping space $\map(BP_D,BSp(2^k)\doscom)_{Bi_{P_D}}$ is homotopy equivalent to $BC_{Sp(2^k)}(P_D)\doscom$,
which is isomorphic to $(B\Z/2)^{2^k}$ (part (1)). Then the obstruction group
\begin{equation*}
H^m(P/P_D;\pi_m(\map(BP_D,BSp(2^k)\doscom))_{Bi_{P_D}})= H^m(P/P_D;\pi_m (B(\Z/2)^{2^k}))
\end{equation*}
is nontrivial only possibly for $m=1$. The group $P/P_D$ is isomorphic
to the group generated by the permutation matrices $B_0,\ldots,B_{k-1}$ (Definition \ref{StubbDef}).
So $P/P_D=(\Z/2)^k$ and the action of $P/P_D$ on $\pi_1 (B(\Z/2)^{2^k})=(\Z/2)^{2^k}$
is given by permutations $\sigma_0,\ldots, \sigma_{k-1}$ which define matrices $B_0,\ldots,B_{k-1}$. 
By Shapiro's lemma \cite[III, Proposition 6.2]{B},
\begin{equation*}
H^1(P/P_D;\pi_m(\map(BP_D,BSp(2^k)\doscom)_{Bi_{P_D}}))=H^1(E_{2^k};(\Z/2)^{2^k})=H^1(1;\Z/2)=1,
\end{equation*}
so all obstruction groups vanish. Therefore $B\alpha$ is homotopic to $Bi_P$ and
$C_{Sp(n)}(\alpha)$ equals $ Z(P)$.

Now we will prove part $(2)$ and $(3)$ for an irreducible $2$-stubborn subgroup $P$ of $Sp(2^k)$.
Let us write $P=Q\wr E_{2^r}$, where $Q$ is an irreducible $2$-stubborn subgroup of $Sp(2^{k-r})$.
Let $\alpha,\beta\colon P\rTo Sp(2^k)$ be two homomorphisms such that $B\alpha^*=
B\beta^*$ and $\alpha|_{BP_T}=i_{P_T}=\beta|_{BP_T}$. We proved that $\alpha|_{BP_D}=i_{P_D}=\beta|_{BP_D}$. 
Let $\bar\alpha,\bar\beta\colon Q^{2^r}\rTo Sp(2^k)$ be the restrictions of $\alpha$ and $\beta$.
Because $Z(Q^{2^r})=Z(Q)^{2^r}=(\Z/2)^{2^r}$, the homomorphisms $\bar\alpha$ and $\bar\beta$ factor
trough homomorpisms 
$$
\widetilde\alpha,\widetilde\beta\colon Q^{2^r}\rTo C_{Sp(2^k)}(Z(Q^{2^r}))=Sp(2^{k-r})^{2^r}.
$$
The map $B\widetilde\alpha$ is homotopic to the map
$$
BQ^{2^r}\simeq \map(BZ(Q^{2^r}),BQ^{2^r})_{Bi}\rTo^{\map(BZ(Q^{2^r}),B\bar\alpha)}
 \map(BZ(Q^{2^r}),BSp(2^k)\doscom)_{Bi}\simeq (BSp(2^{k-r})^{2^r})\doscom,
$$
hence $H^*(B\widetilde\alpha;\F_2)=H^*(\map(BZ(Q^{2^r}),B\bar\alpha);\F_2)$. Analogously
 $H^*(B\widetilde\beta;\F_2)=H^*(\map(BZ(Q^{2^r}),B\bar\beta);\F_2)$.
By Lannes' theory \cite{L}, 
$$
H^*(\map(BZ(Q^{2^r}),B\bar\alpha);\F_2)=
T_{B\bar\alpha^*}^{Z(Q^{2^r})}=T_{B\bar\beta^*}^{Z(Q^{2^r})}=H^*(\map(BZ(Q^{2^r}),B\bar\beta);\F_2),
$$
so $B\widetilde\alpha^*=B\widetilde\beta^*$.

The homomorphisms $\widetilde\alpha$ and $\widetilde\beta$ are matrices of dimension 
$2^r\times 2^r$, where entries are 
$\widetilde\alpha_{i,j},\widetilde\beta_{i,j}\colon Q_i\rTo Sp(2^{k-r})_j$.
The indexes $i$ and $j$ indicate the components in the products.
By induction, $B\widetilde\alpha_{i,i}$ and $B\widetilde\beta_{i,i}$ are homotopic
and therefore $\widetilde\alpha_{i,i}$ and $\widetilde\beta_{i,i}$ are 
conjugate \cite[Th\'eor\`eme 3.4.5]{L}.
We can assume that $\widetilde\alpha_{i,i}=\widetilde\beta_{i,i}$.
Because $Q_i$ and $Q_j$ commutes for $i\ne j$, the homomorphisms 
$\widetilde\alpha_{i,j}$ and $\widetilde\beta_{i,j}$
factor trough homomorphisms 
$\hat\alpha_{i,j},\hat\beta_{i,j}\colon Q_i\rTo C_{Sp(2^{k-r})}(\widetilde\alpha_{j,j}(Q))$.
By induction, the centralizer $C_{Sp(2^{k-r})}(\widetilde\alpha_{j,j}(Q))$ equals $Z(Q_j)=(\Z/2)_j$. 
Because $\widetilde\alpha|_{P_D}=\widetilde\beta|_{P_D}$, the homomorphism 
$\widetilde\alpha_{i,j}\cdot\widetilde\beta_{i,j}^{-1}\colon Q_i\rTo (\Z/2)_j$ factors trough a homomorphism
$\gamma_{i,j}\colon (Q/Q_D)_i\rTo (\Z/2)_j$. Then $\widetilde\beta_{i,j}$ equals to composition
$$
Q_i\rTo^\Delta Q_i\times (Q/Q_D)_i \rTo^{\widetilde\alpha_{i,j}\times\gamma_{i,j}} 
Sp(2^{k-r})_j\times (\Z/2)_j\rTo^\mu Sp(2^{k-r})_j,
$$
where $\Delta$ is the diagonal map composed by the quotient map and $\mu$ is the multiplication in $Sp(2^{k-r})$.
Because $B\widetilde\alpha_{i,j}^*=B\widetilde\beta_{i,j}^*$, the map
$B\gamma_{i,j}$ induces trivial map on mod-$2$ cohomology.
Because $Q/Q_D$ is an iterated wreath product of elementary abelian groups,
the map $\gamma_{i,j}$ is constant \cite[Lemma 6.10]{N}.
Hence $\widetilde\alpha_{i,j}=\widetilde\beta_{i,j}$ and so $\bar\alpha=\bar\beta$ and the 
centralizer $C_{Sp(2^k)}(\alpha)$ is given by the fixed-point set
$C_{Sp(2^k)}(\alpha)=(C_{Sp(2^{k-r})^{2^r}}(\widetilde\alpha))^{E_{2^r}}=
((C_{Sp(2^{k-r})}(Q))^{2^r})^{E_{2^r}}=
((\Z/2)^{2^r})^{E_{2^r}}=\Z/2=Z(P)$, which proves part (2).

The extensions $B\alpha\colon BP\rTo BSp(2^k)\doscom$ of $B\bar\alpha$ are
classified by the obstruction groups $H^m(P/Q^{2^r};\pi_m(\map(BQ^{2^r},BSp(2^k)\doscom)_{B\bar\alpha}))$. By
\cite{JMO}, the mapping space $\map(BQ^{2^r},BSp(2^k)\doscom)_{B\bar\alpha}$ is homotopy equivalent to
$BC_{Sp(2^k)}(Q^{2^r})\doscom=(\Z/2)$. Hence the obstruction groups 
are
$$
H^m(P/Q^{2^r};\pi_m(\map(BQ^{2^r},BSp(2^k)\doscom)_{B\bar\alpha}))=
H^m(P/Q^{2^r};\pi_m(B\Z/2)).
$$ 
The only possible nontrivial obstruction group is for $m=1$.
The group $P/Q^{2^r}=E_{2^r}$ acts by permutation on $BC_{Sp(2^n)}(Q^{2^r})=(B\Z/2)^{2^r}$,
hence $\Ind_{1}^{E_{2^r}} (\Z/2)^{2^r}=(\Z/2)$, by Shapiro's lemma \cite[III, Proposition 6.2]{B}.
Therefore 
$$
H^1(K\wr E_{2^r};(\Z/2)^{2^k})= H^1(K; (\Z/2)^{2^{k-r}}). 
$$ 
Therefore all obstruction groups vanish, so $B\alpha$ is homotopic to $Bi_P$.

Finally let $P=P_1\times\ldots\times P_s$, where $P_i$ is an
irreducible $2$-stubborn subgroup of $Sp(n_i)$.
Let $\alpha,\beta\colon P\rTo Sp(n)$ be two homomorphisms such that
$B\alpha^*=B\beta^*$ and $\alpha|_{P_T}=\beta|_{P_T}$.
Both homomorphisms factor trough $\bar\alpha,\bar\beta\colon P\rTo C_{Sp(n)}(Z(P))=
C_{Sp(n)}((\Z/2)^s)=Sp(n_1)\times\ldots\times Sp(n_s)$. 
In the same way as in the case P is irreducibe $2$-stubborn group,
we can show that $B\bar\alpha^*=B\bar\beta^*$. Maps $\bar\alpha$ and $\bar\beta$
are matrices of dimension $s\times s$ with entries maps 
$\bar\alpha_{i,j},\bar\beta_{i,j}\colon P_i\rTo Sp(n_j)$.
Analogously as before we can show that $\bar\alpha_{i,j}=\bar\beta_{i,j}$, so $B\alpha\simeq B\beta$.
The equation $Z(P)=Z(P_1)\times\ldots\times Z(P_s)$ finishes the proof.
\end{proof}

\section{The map from $Sp(n)\doscom$ to $X$} \label{map}

For every object $Sp(n)/P$ in $\widetilde{\mathcal R}(Sp(n))$ we define a map
$f_P\colon P\rTo X$ as the composition of the two inclusions $i_P\colon P\rTo N$ and
$f_N\colon N\rTo X$. We will prove that for every morphism $c_g\colon Sp(n)/P\rTo Sp(n)/Q$ in
$\widetilde{\mathcal R}(Sp(n))$, the diagram
\begin{align} \label{PokaziKom}
\begin{diagram}
 BP & \rTo^{Bi_P} & BN &\rTo^{Bf_N}& BX \\
 \dTo<{Bc_g} &&&& \dIgual \\
 BQ & \rTo^{Bi_Q} &BN &\rTo^{Bf_N}& BX \\
\end{diagram}
\end{align}
commutes up to homotopy.

Let us define $\alpha=f_N\circ i_P$ and $\beta=f_N\circ i_Q\circ c_g$.
Then $B\alpha^*=B\alpha^*$.  The group $P_T=P\cap
Sp(n)$ is $2$-toral.  The restrictions $\alpha|_{P_T}$ and
$\beta|_{P_T}$ are conjugate in $Sp(n)$, and hence by
\cite[Proposition 4.1]{MN1}, they are also conjugate in the normalizer
$N$ of the maximal torus.  So $B\alpha|_{BP_T}\simeq B\beta|_{BP_T}$. By the
next proposition, $B\alpha\simeq B\beta$.

Let $K\rTo G\rTo H$ be an exact sequence of groups. Then $H$ acts freely on
$\widetilde{BK}=EG/K\simeq BK$, and $\widetilde{BK}/H$ equals $BG$. For any space
$BX$ with trivial action of the group $H$, we have
\begin{multline} \label{hH}
\map(BG,BX)=\map(\widetilde{BK}/H,BX)=\map_H(\widetilde{BK},BX)\simeq\\
\simeq\map_H(EH\times \widetilde{BK},BX)=\map_H(EH,\map(\widetilde{BK},BX))=
\map(\widetilde{BK},BX)^{hH}.
\end{multline}

\begin{Prop}
For every $Sp(n)/P\in ob(\widetilde{\mathcal R}(Sp(n)))$, the canonical map
$$
\pi_0(\map(BP,BX)_{B\alpha|_{BP_T}=Bf_{P_T}})\rTo Hom(H^*BX,H^*BP)
$$
is an injection.
\end{Prop}
\begin{proof}
Consider the diagram
\begin{diagram}\label{prvi}
&& \map(\widetilde{BP_T},BY)_{Bi_T}&&\\ & \ldTo&&\rdTo&\\
\map(\widetilde{BP_T},BSp(n)\doscom)_{Bi_T}&&&&\map(\widetilde{BP_T},BX)_{Bf_T}.
\end{diagram}
By \cite{N5}, the mapping space $\map(\widetilde{BP_T},BSp(n)\doscom)_{Bi_T}$
is homotopy equivalent to $BC_{Sp(n)}(BP_T)\doscom$ and by proposition \ref{5.5},
this space is homotopy equivalent to $(BT_{Sp(n)})\doscom$. Analogously
$\map(\widetilde{BP_T},BY)_{Bi_T}$ is homotopy equivalent to $(BT_{Sp(n)})\doscom$.
The mapping space $\map(\widetilde{BP_T},BX)_{Bf_T}$ is the classifying space
of a $2$-compact group (\cite{DW}). Its Weyl group is
$Iso(Bf_N\circ Bi_{P_T}t)=\{w\in W_X\mid w\circ Bf_N\circ Bi_{P_T}\simeq Bf_N\circ Bi_{P_T}\}$
(\cite[Proposition 4.3]{V}). By the construction of the map $f_N$, the group $Iso(Bf_N\circ Bi_{P_T})$
equals $Iso(Bi_N\circ Bi_{P_T})$. Because $Iso(Bi_N\circ Bi_{P_T})$ is the Weyl group of
the mapping space $\map(\widetilde{BP_T},BSp(n)\doscom)_{Bi_T}\simeq (BT_{Sp(n)})\doscom$, 
the group $Iso(Bf_N\circ Bi_{P_T})$ is trivial, hence $\map(\widetilde{BP_T},BX)_{Bf_T}
\simeq (BT_{Sp(n)})\doscom$.
Therefore both maps in the diagram (\ref{prvi}) are homotopy equivalences.

Taking homotopy fixed points we obtain the following diagram
\begin{diagram}
&&\map(\widetilde{BP_T},BY)_{Bi_T}^{h(P/P_T)}&&\\ &\ldTo&&\rdTo&\\
\map(\widetilde{BP_T},BSp(n)\doscom)_{Bi_T}^{h(P/P_T)}&&&&\map(\widetilde{BP_T},BX)_{Bf_T}^{h(P/P_T)},
\end{diagram}
where both maps are mod-$2$ equivalences, since an equivariant mod-$2$ equivalence
between 1-connected spaces induces a mod-$2$ equivalence between the homotopy
fixed-point sets. 

By proposition \ref{5.5}, the components of $\map(\widetilde{BP_T},BSp(n)\doscom)_{Bi_T}^{h(P/P_T)}$
are distinguished by mod $2$ cohomology. Any map in 
$\map(\widetilde{BP_T},BX)_{Bf_T}^{h(P/P_T)}$ has a lift to $BN$ and therefore to $BY$.
The obstruction group which classifies the extensions is
$$
H^2(P/P_T;\pi_2\map(\widetilde{BP_T},BX)_{Bf_T})\cong
H^2(P/P_T;\pi_2\map(\widetilde{BP_T},BSp(n)\doscom)_{Bi_T}),
$$
so the components of $\map(\widetilde{BP_T},BX)_{Bf_T}^{h(P/P_T)}\simeq 
\map(BP,BX)_{B\alpha|_{BP_T}= Bi_{P_T}}$ are also distinguished by mod-$2$ cohomology.
\end{proof}


Diagram \eqref{PokaziKom} establishes a map from the $1$-skeleton of the homotopy colimit
$\{BP\}_{\widetilde{\mathcal{R}}_2(Sp(n))}$ to $BX$. The obstruction groups for extending a map
defined on the $1$-skeleton of the homotopy colimit to a map on the total homotopy colimit are $$
\sideset{}{^{i+1}}\lim_{\overleftarrow{\widetilde{\mathcal{R}}_2(Sp(n))}}
\pi_i\map(BP,BX)_{Bf_P}$$ for $i\ge 2$, where $\varprojlim^i$ is the $i$-th derived functor of the
inverse limit functor (\cite{BK} and \cite{WO}).

Let $\mathcal{A}b$ be the category of Abelian groups and let $$ \Pi_j^X,\Pi_j^{Sp(n)}\colon
\widetilde{\mathcal{R}}_2(Sp(n)) \rTo\mathcal{A}b $$ be functors defined by $$
\Pi_j^X(Sp(n)/P)=\pi_j\map(BP,BX)_{Bf_P},$$ $$ \Pi_j^{Sp(n)}(Sp(n)/P)=
\pi_j\map(BP,BSp(n)\doscom)_{Bi_P}.$$ Note that $\map(BP,BSp(n)\doscom)_{Bi_P}$ is isomorphic 
to $BZ(P)\doscom$ (\cite{JMO}, Theorem 3.2) therefore $\Pi_1(Sp(n))(Sp(n)/P)$ is well defined. 
By the next proposition, also $\Pi_1(X)(Sp(n)/P)$ is well defined.

\begin{Prop}
There exists a natural transformation $$ \mathcal{T}\colon \Pi_j^{Sp(n)}\to\Pi_j^X $$ which is an
equivalence.
\end{Prop}
\begin{proof}
For every $2$-stubborn group $P$ we have homotopy
equivalences
\begin{align} \label{zvezdica}
\map(BP,BSp(n)\doscom)_{Bi_P}\lTo^{\simeq}\map(BP,BY)_{Bi_P}
\rTo^{\simeq}\map(BP,BX)_{Bf_P}
\end{align}
which depend on the chosen lift $Bi_P\colon BP\rTo BY$ of the map $Bi_P\colon BP\rTo
BSp(n)\doscom$. Because $\Rep(P,Sp(n)) \rTo \left[ BP,BSp(n) \right] $ is a bijection (\cite{DZ},
\cite{N5}), two lifts differ by a conjugation $Bc_g$. Since $Bf_p\simeq Bf_P\circ Bc_g$, the
equivalence (\ref{zvezdica}) induces well defined isomorphisms $$
\Pi_j^{Sp(n)}(Sp(n)/P)\to\Pi_j^X(Sp(n)/P) $$ which commute with maps induced by morphisms in
$\widetilde{\mathcal{R}}_2(Sp(n))$.
\end{proof}

\begin{Prop}
For all $i,j\ge 1$,
$$\sideset{}{^i}\lim_{\overleftarrow{\widetilde{\mathcal{R}}_2(Sp(n))}}
\pi_j\map(BP,BX)_{Bf_P}=0.$$
\end{Prop}
\begin{proof}
By the previous lemma,
$$\sideset{}{^i}\lim_{\overleftarrow{\widetilde{\mathcal{R}}_2(Sp(n))}}
\pi_j\map(BP,BX)_{Bf_P}= \sideset{}{^i}\lim_{\overleftarrow{\widetilde{\mathcal{R}}_2(Sp(n))}}
\pi_j\map(BP,BSp(n)\doscom)_{(Bi_P)\doscom}$$ and the right side is $0$ \cite[Theorem 4.8]{JMO}.
\end{proof}

Because all obstructions vanish, there exists a map $$ f\colon \hocolim_{{\widetilde
R}_2(Sp(n))} BP\rTo BX. $$ By the construction of the map we have a commutative
diagram
\begin{diagram}
&&BY &&\\ &\ldTo&&\rdTo& \\ \hocolim_{{\widetilde R}_2(Sp(n))} ESp(n)/P&&\rTo^{f} &&BX
\end{diagram}
where the diagonal maps induce monomorphisms on the cohomology and therefore $f^*$ is a
monomorphism. Since $H^*BSp(n)\cong H^*BX$, $f^*$ is an isomorphism and therefore $f$ is a
homotopy equivalence.

\section{$Sp(n)$ as a loop space} \label{loop}

The normalizer conjecture can be stated also for finite loop spaces with maximal torus
normalizers as a weak version of Wilkerson' conjecture (see \cite{W}).

\begin{Thm}
Let $L$ be a connected finite loop space with a maximal torus normalizer isomorphic to
that of $Sp(n)$. Then $BL$ is homotopy equivalent to $BSp(n)$.
\end{Thm}
\begin{proof}
To prove $BL\simeq BSp(n)$ is equivalent to showing that $BL$ and $BSp(n)$ lie in the
same adic genus (\cite{N6}). The loop spaces $BL$ and $BSp(n)$ have the same rational
genus. Since $BL$ is finite and connected, $L\pcom$ is a $p$-compact group. The
maximal torus normalizer of $L\pcom$ is just the fibrewise $p$-completion of $N$ by
the fibration $BT\rTo BN\rTo BW_{L}$. Hence $L\pcom$ and $Sp(n)\pcom$ have isomorphic
normalizers of the maximal torus. By \cite{AGMV}, $BSp(n)\pcom$ is $N$-determined if
$p$ is an odd prime and by the main theorem of this paper $BSp(n)\doscom$ is (weakly)
$N$-determined. So $BL\pcom$ and $BSp(n)\pcom$ are homotopy equivalent.
\end{proof}



\begin{thebibliography}{ABCD}

\bibitem{A} J. Aguad\'e, {\it Constructing modular classifying spaces}, Israel J. Math.
{\bf 66} (1989), 23--40.

\bibitem{AGMV} K. Andersen, J. Grodal, J.M. M\o ller, A. Viruel, {\it The classification
of $p$-compact groups, $p$ odd}, Preprint http://arxiv.org/abs/math.AT/0302346.


\bibitem{B} K.S. Brown, {\it Cohomology of groups} Springer-Verlag New
York Berlin Heidelberg London Paris Tokyo Hong Kong Barcelona Budapest (1994).

\bibitem{BK}  A. Bousfield, D. Kan, {\it Homotopy limits, completion and
localisation}, SLNM {\bf 304}, Springer Verlag (1972).

\bibitem{CE}  A. Clark, J. R. Ewing, {\it The realization of polynomial algebras
as cohomology rings}, Pacific J. Math. {\bf 50} (1974), 425--434.


\bibitem{CWW} M.\ Curtis, A.\ Wiederhold, B.\ Williams, {\it Normalizers of maximal tori}, Lecture
Notes in Math.\ {\bf 418}, 31--47.

\bibitem{D}, W.G. Dwyer, {\it Lie groups and $p$-compact groups},
Proceedings of the International Congress of Mathematicians, Vol. II (Berlin, 1998).
Doc. Math. 1998, Extra Vol. II, 433--442 (electronic).


\bibitem{DW1} W.G. Dwyer, C.W. Wilkerson, {\it The center of a $p$-compact group},
The \v Cech centennial (Boston, MA, 1993), 119--157, Contemp. Math., 181, Amer. Math.
Soc., Providence, RI, 1995.


\bibitem{DW} W.G. Dwyer, C.W. Wilkerson, {\it Homotopy fixed point
methods for Lie groups and finite loop space}, Ann. Math. {\bf 139} (1994), 395--442.

\bibitem{DW2} W.G. Dwyer, C.W. Wilkerson, {\it $p$-compact groups with abelian Weyl groups}, Preprint.

\bibitem{DW3} W.G. Dwyer, C.W. Wilkerson, {\it A new finite loop space at the prime two}, J. Amer. Math. Soc.
{\bf 6} (1993), 37--63.

\bibitem{DW5} W.G. Dwyer, C.W. Wilkerson, {\it Product splitting for $p$-compact groups},
Fund. Math. {\bf 147} (1995), no. 3, 279--300.

\bibitem{DW6} W.G. Dwyer, C.W. Wilkerson, {\it A cohomology decomposition theorem}, Topology
{\bf 31} (1992), 433--443.

\bibitem{DZ} W.G. Dwyer, A. Zabrodsky, {\it Maps between classifying spaces},
Algebraic topology, Barcelona, 1986, 106--119, Lecture Notes in Math., 1298, Springer,
Berlin, 1987.

\bibitem{GLZ} J.H.\ Gunawardena, J.\ Lannes, S.\ Zarati, {\it Cohomology des groupes sym\' etrques
et application de Quillen}, Advances in homotopy theory (Cortona, 1988) 61-68, London
Math. Soc. Lecture Note Ser., 139.

\bibitem{JMO} S.\ Jackowski, J. McClure, R. Oliver, {\it
Homotopy classification of self-maps of $BG$ via $G$-actions, part I and part II},
Ann. Math. {\bf 135} (1992), 183--270.


\bibitem{L} J.\ Lannes, {\it Sur les espaces fonctionelles dont la source
est la classifiant d'un $p$-groupe abe\'elien \'elem\'entaire} Publ. Math. IHES {\bf
75} (1992), 135--244.


\bibitem{M} J.M.\ M\o ller, {\it Normalizers of maximal tori} Math. Z. {\bf 231}
(1999), no. 1, 51--74.

\bibitem{M1} J.M.\ M\o ller, {\it Racional isomorphism of $p$-compact groups},
Topology {\bf 35} (1996), no. 1, 201--225.

\bibitem{M2} J.M.\ M\o ller, {\it Deterministic $p$-compact groups},
Stable and unstable homotopy (Toronto, ON, 1996), 255--278, Fields Inst.\ Commun., 19, Amer. Math.
Soc., Providence, RI, 1998.

\bibitem{M3} J.M. M\o ller, {\it Homotopy Lie groups}, Bull. Amer. Math. Soc.
{\bf 32} (1995), no. 4, 413--428.

\bibitem{MN} J.M. M\o ller, D. Notbohm, {\it Centers and finite
coverings of finite loop spaces}, J. Reine Angew. Math. 456 (1994), 99--133.

\bibitem{MN1} J.M. M\o ller, D. Notbohm, {\it Connected finite loop space with maximal tori},
Trans. Amer. Math. Soc. {\bf 350} (1998), no. 9, 3483--3504.


\bibitem{N} D. Notbohm, {\it Homotopy uniqueness of classifying spaces
of compact connected Lie groups at primes dividing the order of the Weyl group},
Topology {\bf 33} (1994), no. 2, 271--330.

\bibitem{N1} D. Notbohm, {\it A uniqueness result for orthogonal groups as
$2$-compact groups}, Arch. Math. (Basel) {\bf 78} (2002), no. 2, 110--119.

\bibitem{N55} D. Notbohm, {\it On the $2$-compact group $DI(4)$}, Preprint.


\bibitem{N3} D. Notbohm, {\it Spaces with polynomial mod-$p$ cohomology}, Proc. Cambridge Philos. Soc.
{\bf 126} (1999), 277--292.

\bibitem{N4} D. Notbohm, {\it Classifying spaces of compact Lie groups and finite loop spaces},
Handbook of Algebraic Topology, North-Holland (1995).

\bibitem{N5} D. Notbohm, {\it Maps between classifying spaces}, Math. Z. {\bf 207} (1991), 229--257.

\bibitem{N6} D. Notbohm, {\it Fake Lie groups with maximal tori IV},
Math. Ann. {\bf 294} (1992), no. 1, 109--116.

\bibitem{O} R. Oliver, {\it $p$-stubborn groups of the classical compact Lie
groups}, journal of Pure and Applied Algebra {\bf 92}, 55--78.

\bibitem{Q} D. Quillen, {\it On the cohomology and $K$-theory of the general linear group over a finite
field}, Ann. Math. {\bf 96} (1972), 552--586.

\bibitem{VV} A. Vavpeti\v c, A. Viruel, {\it On the homotopy type of the classifying space
of the exeptional Lie group of rank 4}, Manuscripta Math. {\bf 107} (2002), no. 4,
521--540.

\bibitem{V} A. Viruel, {\it Homotopy uniquness of $BG_2$}, Manuscripta Math. {\bf 95}
(1998), no. 4, 471--497.

\bibitem{W} C.W. Wilkerson, {\it Rational maximal tori},
J. Pure Appl. Algebra {\bf 4} (1974), 261--272.

\bibitem{WO} Z. Wojtkowiak, {\it On maps from holim $F$ to $Z$}, in
Algebraic Topology, Barcelona 1986, SLNM {\bf 1298}, 227--236.

\bibitem{Z} A. Zabrodsky, {\it On the realization of invariant subgroups of $\pi_*(X)$},
Trans. Amer. Math. Soc. {\bf 285} (1984), 467--496.

\end{thebibliography}
\end{document}